\documentclass[numbook, envcountsect, envcountsame, envcountreset, runningheads, smallextended]{svjour3}

\usepackage{amscd,amsfonts,amssymb,amsmath,latexsym,array,hhline,xcolor,graphicx}

\usepackage{latexsym}
\usepackage{times}

\newcommand\cE{{\cal E}}
\newcommand\cC{{\cal C}}
\newcommand\cF{{\cal F}}

\newcommand\cM{{\cal M}}
\newcommand\cX{{\cal X}}

\newcommand\cZ{{\cal Z}}
\newcommand\e{{\varepsilon}}

\def\E{{\bf E}}

\def\text#1{\hbox{#1}}

\def\E{{\bf E}}

\def\build #1_#2{\mathrel{\mathop{\kern 0pt #1}\limits_{#2}}} 
\newcommand{\zs}[1]{{\mathchoice{#1}{#1}{\lower.25ex\hbox{$\scriptstyle#1$}}
{\lower0.25ex\hbox{$\scriptscriptstyle#1$}}}}

\numberwithin{equation}{section}

\newcommand\fdem{$\Box$}
\newcommand\beq{\begin{equation}}
\newcommand\eeq{\end{equation}}
\newcommand\bea{\begin{eqnarray}}
\newcommand\eea{\end{eqnarray}}
\newcommand\bean{\begin{eqnarray*}}
\newcommand\eean{\end{eqnarray*}}

\newcommand\beal{\begin{alighn}}
\newcommand\eeal{\end{align}}

\begin{document}

\title{ On the entropy-minimal  martingale measure in the exponential Ornstein--Uhlenbeck stochastic volatility model}
	\author{ Yuri Kabanov \and  Mikhail A. Sonin} 
	
	\institute{\at	Lomonosov Moscow State University, Vega Institute Foundation, 
	Moscow, Russia, and \\
	 Universit\'e de Franche-Comt\'e, Laboratoire de Math\'ematiques, UMR CNRS 6623, 
	Besan\c{c}on,  France \\
  \email{ykabanov@univ-fcomte.fr}. \\
	\and Lomonosov Moscow State University, Moscow, Russia\\
	 \email{msonin1999@yandex.ru}}

\titlerunning{On the entropy-minimal  martingale measure}
\date{\today}
\maketitle

\qquad \qquad \qquad \qquad \qquad \qquad 
{\sl Dedicated to the 80th anniversary  of Isaak Sonin.}\\
\\

\begin{abstract}
We consider a  stochastic volatility model  where  the price evolution  depend on the exponential of the Ornstein--Uhlenbeck process. After a  brief revision of the related theory  the entropy-minimal equivalent martingale measure. 
is calculated. 
\end{abstract} 

\keywords{
option pricing \and  stochastic volatility \and  Ornstein--Uhlenbeck process \and entropy-minimal martingale measure \ 
}
\smallskip

\noindent
 {\bf Mathematics Subject Classification (2010)} 60G44

\section{Introduction}  The modern theory of stochastic processes is highly indebted to quantitative finance which is  the vast field of its applications. 
Unlike the physical world governed by the eternal laws of the nature, ``laws" of  finance permanently 
vary due to various factors including even a theoretical development which may strongly influence the behavior of economic agents.  There is a non-stop competition of models describing the dynamics of basic securities as well as ideas how to price the derivative securities. It may happen  that specific features  that was considered yesterday as  drawbacks turn to be advantages today as reflecting new reality.   In their everyday practice, financial engineers (``quants") are asked questions what is the price of such and such option contract.  

Mathematics provided a plethora of models and  general pricing recommendations.  The simplest one is the classic  "pricing by replication" paradigm  suggested 
by Black and  Scholes for a model where price of the basic securities are given by  geometric Brownian motions and the number of ``sources of randomness" is equal to the number  of risky assets. In this case the ``market", or better to say, market model is complete, the equivalent martingale measure is unique,  and every European-type contingent claim can be replicated by the terminal value a portfolio. 
Unfortunately,  in the majority of cases statistical tests reject the hypothesis that the logarithms of price increments  are Gaussian random variables.  This fact leads to a search for better  models  of the price dynamics.  
 
 Amongst  models of price evolution of basic securities, nowadays the so-called stochastic volatility models again become   very popular. One of them, sometimes referred to as the exponential Ornstein--Uhlenbeck model, see Mejia Vega \cite{Vega},  is considered in the present note. 
 
 The general recommendation how to provide an arbitrage-free price of an option consists in computation of its mean with respect to an equivalent martingale measure (EMM). The delicate issue here is with respect to which one. 
 Various solutions are suggested. Some experts recommend to enlarge the framework  by considering more traded securities 
 (for example, traded ``vanilla" options) to get a complete market with  a single equivalent martingale  measure. 
 In other words, there is a believe that market follows the "pricing by replication"  rule.  

Other  experts defend the idea that 
financial  markets choose a relevant  equivalent martingale measure using a certain variational principle. A reasonable choice is the entropy minimization one.  
We compute such a martingale  measure for the considered model and use this occasion to present  a brief, hopefully, comprehensive  piece  of  the beautiful theory of entropy minimal martingale measures usually suffered from technical details and generalizations. In our presentation we are based mainly on papers  by Hobson \cite{Hobson}, Frittelli \cite{Frittelli}, and  Kabanov and Stricker, \cite{KS}. In our numerical experiments we observe 
using real data (but artificially chosen coefficients) that  such a pricing rule gives rather good approximation of the observed option prices.

\section{Model}
We consider the option pricing for two asset  stochastic volatility model where the riskless asset is taken as the numeraire and  $S=(S_t)$, the  price process of a risky asset, is given 
by the system of SDEs
\bea
\label{1}
dS_t&=&S_t(\mu(t,V_t)dt+ \sigma(V_t)dB_t)\\
\label{2}
dV_t&=&(\theta - \alpha V_t)dt+ \beta dW_t, \quad V_0=v,
\eea 
where $W=\rho B+\bar \rho \tilde B$,   $B$ and $\tilde B$ are independent Wiener processes,   
$\rho\in [-1,1]$, $\bar \rho:=\sqrt {1-\rho^2}$. 
The coefficients have a specific form: 
 \beq
 \label{3}
 \mu(t,V_t)=(\mu+kt)e^{V_t}, \qquad  \sigma(V_t)=\sigma e^{V_t}. 
 \eeq
The parameters  $\mu$, $\kappa$, and $\theta$ are arbitrary real numbers, $\sigma, \alpha, \beta>0$.  

The mean reverting Ornstein--Uhlenbeck process $V$ is Gaussian. It admits the explicit formula   
\beq
V_t=v e^{-\alpha t}+  \frac{\theta}{\alpha} \big(1-e^{-\alpha t}\big)+ \beta \int_0^te^{-\alpha (t-s)} dW_s  
\eeq
implying that $
V_t\sim {\cal N}\left(\theta/\alpha + e^{-\alpha t}(v - {\theta}/{\alpha}), \beta^2/(2\alpha)(1 - e^{-2\alpha t})\right)$. 

The process $e^V$ sometimes is called exponential Ornstein--Uhlenbeck process. 
 As the usual Ornstein--Uhlenbeck process, it has the mean-reverting property. There are statistical evidences showing that this property is inherent to the volatility of financial assets, see Fig. 8.1. Since volatility is positive, it seems  reasonable to use the model  $\sigma (V_t) = \sigma e^{V_t}$
where the scale parameter $\sigma>0$ adds a flexibility.

The model has two sources of randomness, the equivalent martingale measure is not unique and one   cannot represent every contingent claim by a stochastic integral with respect to the price process $S$.  The Back--Scholes paradigm of ``pricing by replication" fails.  There is no consensus between experts how to price options is such a situation. However, there is a believe that the market prices of  options can be obtained by calculating  expectations   
of  corresponding pay-offs with respect to an equivalent martingale measure satisfying some variational property. 
The most popular choices are equivalent entropy minimal martingale measures and variance minimal martingale  
measures. 

\section{Basic definitions and properties}
Let $(\Omega, {\bf F}, (\cF_t)_{t\le T}, P)$ be a filtered probability space and let $\cX=\{X\}$ be a family  of continuous semimartingales.  We define the set  $\cM$  (resp., $\cM^e$) of absolute continuous martingale measures (resp. 
equivalent martingale measures) as the 
set of probability measures   
 $Q\ll P$ (resp., $Q\sim P$)  such  that  each  $X\in \cX$ is a  martingale with respect to $Q$.  
 We denote $\cZ$ and $\cZ^e$ the corresponding sets of Radon--Nikodym derivatives $Z=dQ/dP$. Recall that 
  the (relative) entropy of $Q$ with respect to $P$ is $I(Q,P):=\E Z\ln Z =\E^Q \ln Z $. 
    
  The {\sl entropy minimal equivalent martingale measure} $Q^o$ (EMM) is defined as an equivalent martingale measure  such that 
 $$
 I(Q^o,P)=\inf_{Q\in \cM^e} I(Q,P).
 $$ 
  It is more convenient to define $EMM$ as an element $Q^o\in \cM^e$ with the density $Z^o\in \cZ^\e$ such that 
$$
\E Z^o\ln Z^o=J^o:=\inf_{Z\in \cZ^e}\E Z\ln Z.
$$ 

\section{Prerequisites from the general theory}

Let $\cZ$ be a non-empty convex set of random variables   $Z\ge 0$ with  $\E Z=1$.  We introduce two i convex subsets of $\cZ$, 
$$
\cZ^e:=\{Z\in \cZ\colon Z>0\}, \qquad \cZ_\varphi:=\{Z\in \cZ\colon \E \varphi(Z)<\infty\},
$$ 
where 
$\varphi(x)=x\ln x$, $x> 0$, $\varphi(0)=0$. The function $\varphi(x) \ge -e^{-1}$, it  is  convex,  its derivative $\varphi'(x)=\ln x+ 1$, 
$\varphi'(0)=-\infty$, and   $\psi(x):=e^{x-1}$ is the inverse of $\varphi'(x)$. 

\smallskip

For  $Z,\tilde Z\in \cZ_\varphi$ and $t\in [0,T]$  we put $f_t:=\varphi (t Z+(1-t)\tilde Z)$ and $F_t:=\E f_t$. 

\smallskip

Note that $\E Z\ln Z<\infty$ if and only if $\E Z\ln^+ Z<\infty$. Moreover, 
$$
\E Z\ln^+ \tilde Z=\E I_{\{Z\ge \tilde Z\}}Z\ln^+ \tilde Z+ \E I_{\{Z< \tilde Z\}} Z\ln^+ \tilde Z\le  \E \varphi(Z)+\E \varphi(\tilde Z)+2e^{-1}<\infty. 
$$
Similarly, $\E Z\ln^+ \tilde Z\le  \E \varphi(Z)+\E \varphi(\tilde Z)<\infty $.  Thus the random  variable $Z\ln \tilde Z$ in this framework  is  integrable.  

\begin{lemma} $F'_0=\E f'_0=\E Z\ln \tilde Z-\E\varphi (\tilde Z)$. 
\end{lemma}
{\sl Proof.} 
Let  $0<s\le t\le 1$. Then  
$(f_t-f_0)/t\ge (f_s-f_0)/s$. Indeed, this inequality  is equivalent to the inequality $f_s\le (s/t)f_t+(1-s/t)f_0$ which holds since $f$ is convex. Also $f_1-f_0=\varphi (Z)-\varphi (\tilde Z)$ where the rhs is an integrable r.v. By the monotone convergence  $F'_0=\E f'_0$. It remains to note that  
$$
f'_0:=\varphi' (\tilde Z)(Z-\tilde Z)=Z-\tilde Z+Z\ln \tilde Z-\varphi (\tilde Z)  
$$
and, therefore, $F'_0=\E Z\ln \tilde Z-\E\varphi (\tilde Z)$. 
\fdem 

\smallskip
Put $\widehat J^o:=\inf_{Z\in \cZ}\E \varphi(Z)$. It happens that if  a minimizer exists (necessarily unique because of the strict convexity of $\varphi $), then 
it is an element of  $\cZ^e$ provided that the latter set  is non-empty. The formal statement is 

\begin{lemma} 
If $\widehat J^o=\E \varphi(Z^o)$ for some $Z^o\in \cZ_\varphi$ and  $\cZ^e\cap \cZ_\varphi\neq \emptyset$, then $Z^o\in \cZ^e$ and $\widehat J^o=J^o$.  
\end{lemma}
{\sl Proof.} Let $Z\in \cZ^e\cap \cZ_\varphi$.  Applying the above lemma (with $\tilde Z=Z^o$)  we get that 
$F'_0=\E Z\ln Z^o-\widehat J^o$. Since the convex function $F_t$ attains its minimum at zero, $F'_0\ge 0$.  But this is possible only if $Z^ o >0$ (a.s.). \fdem

\smallskip

 Let $
 \cC:=\{\eta\colon \E |\eta|Z<\infty\ \hbox{and\ } \E\eta Z \le 0\ \forall\, Z\in \cZ_\varphi\}$.

 \smallskip
 The next criteria provides a checking tool.  
 
 \begin{lemma} 
 \label{A3}
Let $Z^ o\in  \cZ_\varphi\cap \cZ^e$ and let $J:=\E \varphi(Z^o)$. Then  
$$
J^o=J\Longleftrightarrow  \hbox{ there is $\eta\in \cC$
such that $\E\eta Z^o=0$ and $Z^o=e^{J-\eta}$}. 
$$
\end{lemma}
{\sl Proof.} ($\Leftarrow$) Due to convexity of $\varphi $ we have for any $Z\in \cZ_\varphi$ that 
$$
\E\varphi (Z)\ge \E\varphi (Z^ o)+\E\varphi' (Z^o)(Z-Z^o)=\E\varphi (Z^ o)+\E(1+J-\eta)(Z-Z^o)\ge \E\varphi (Z^ o). 
$$

($\Rightarrow$) Let $J^o=\E \varphi(Z^o)$. Let $\eta=J^o-\ln Z^o$. Then $\E\eta Z^o=0$ and for any $Z\in \cZ_\varphi$
$$
\E \eta Z=J^o-\E Z\ln Z^o=\E \varphi (Z^o)-\E Z\ln Z^o=-F'_0\le 0 
$$
since $Z^o$ is the minimizer. \fdem

\begin{lemma} 
If $\cZ$ is closed in $L^0$ and $\widehat J^o<\infty$, then there is $Z^o\in \cZ$ such that 
$J^o=\E \varphi(Z^o)$. 
\end{lemma}

\smallskip 
The above formulation is slightly more general than Prop. 3.2 in \cite{KS}, but it could be obtained 
by the similar arguments where the reference to the Koml\'os theorem is replaced by the reference to the  
von Weizs\"acker theorem, see \cite{KabS}, Th. 5.2.3. 

\smallskip
\noindent
{\sl Remark. } The presentation of this section is a simplified version of that in \cite{KS} where more general functions than $x\ln x$ are considered.

\section {The Hobson construction}
Let us consider the  system of SDE, more general than (\ref{1}) - (\ref{2}): 
\bea
\label{1+}
dS_t&=&S_t Y_t(\lambda_tdt+ dB_t),\\
\label{2+}
dY_t&=&\alpha_tdt+ \beta_t dW_t, \quad V_0=v,
\eea 
where $\lambda_t = \lambda(t,Y_t)$, $\alpha_t = \alpha(t,Y_t)$, $\beta_t = \beta(t,Y_t)$ where $\lambda,\alpha, \beta$ are continuous functions. We assume that 
 the second equation above has a unique strong solution and $\lambda^2\cdot\Lambda_ T<\infty$.

Let $\xi$ be a predictable process with $\xi^2\cdot \Lambda_T<\infty$. 
Define a strictly positive local martingale 
$$
Z^\xi =\cE(-\lambda \cdot B-\xi \cdot \tilde B)
=\exp \left\{-\lambda\cdot B -\frac 12 \lambda^2\cdot \Lambda  -\xi\cdot \tilde B-\frac 12 \xi^2\cdot \Lambda\right\}.
$$
Here and the sequel  we use, where it is convenient, the standard notations for integrals (ordinary and stochastic); $\Lambda_t\equiv t$. 

\smallskip
If $\E Z^\xi_T=1$, then  $Q^\xi=Z^\xi_T P$ is a probability measure equivalent to $P$. According to the Girsanov teorem, under $Q^\xi$  both
$B^\lambda:=B+\lambda \cdot \Lambda$ and $\tilde B^\xi:=\tilde B+\xi\cdot \Lambda$ are Wiener processes.  

\smallskip

The following result is a part of Th. 3.1 from \cite{Hobson}. 

\begin{proposition} 
\label{EMM}
Let  $H$ be a predictable process with $H^2\cdot \Lambda_T<\infty$ be such that 
\beq
\label{Hequation}
\frac 12 \lambda^2\cdot \Lambda_T-H \cdot B^\lambda_T-\xi\cdot \tilde B_T-\frac 12 \xi^2\cdot \Lambda_T =  J={\rm const}. 
 \eeq

$(a)$ Suppose that $(H-\lambda)\cdot B^{\lambda}$ is a $Q^\xi$-martingale.  Then $I(Q^\xi,P)=J$. 

$(b)$ If, moreover,  the process
$(H-\lambda)\cdot B^\lambda$ is a $Q$-martingale with respect to any equivalent 
martingale measure $Q$ with finite entropy with respect to $P$, 
then $Q^\xi$ is the entropy minimal EMM and  $I(Q^\xi,P)=J=J^o$. 
\end{proposition} 
{\sl Proof.} $(a)$ Let  $\E^\xi$ denote the expectation with respect to the measure $Q^\xi$. Then $\E^\xi(H-\lambda)\cdot B^\lambda_T=0$ and we get after regrouping terms that 
\bean
\E Z^\xi_T\ln Z^\xi_T&=& \E^\xi \ln Z^\xi_T=\E^\xi \Big[-\lambda \cdot B_T-\frac 12 \lambda ^2\cdot\Lambda _T
-\xi\cdot \tilde B_T-\frac 12 \xi^2\cdot \Lambda_T\Big]\\
&=&\E^\xi \Big[(H-\lambda) \cdot B^\lambda_T+\frac 12 \lambda ^2\cdot\Lambda _T-H \cdot B^\lambda_T
-\xi\cdot \tilde B_T-\frac 12 \xi^2\cdot \Lambda_T\Big]=J.
\eean

$(b)$ According to (\ref{Hequation}), $Z^\xi_T=e^{J-\eta}$ where
$\eta:=(\lambda-H)\cdot B^\lambda_T$. Due to our assumption $\E^Q\eta=0$ for any equivalent martingale measure with finite entropy and the result follows from Lemma \ref{A3}. \fdem 

\section{Solution of the Hobson equation}

    In the Hobson equation (\ref{Hequation}) we should determine the processes $\xi$ and $H$ and the constant $J$. 
To this end,  we consider  on $[0,T]$ the  terminal value problem  
\beq
\label{PDE}
\dot f+\frac 12 \beta^2 f''-\frac 12 \bar \rho^2\beta^2f'^2+(\alpha-\rho\beta\lambda)f'+\frac 12\lambda^2=0, \qquad f(T,y)=0,  
\eeq
where the functions $\alpha,\beta,\lambda, \rho, \bar\rho $ of the variable $(t,y)$ are known, $\dot f$ and $f'$ are derivatives in $t$ and $y$, respectively. 
\begin{proposition} Suppose that the problem (\ref{PDE}) has a solution $f\in C^{1,2}$. 
Then 
$$
\xi_t= \rho(t,Y_t)\beta (t,Y_t) f'(t,Y_t), \qquad H_t=\bar \rho(t,Y_t)\beta (t,Y_t) f'(t,Y_t) ,\quad  J=f(0,Y_0). 
$$
\end{proposition}
{\sl Proof.} Recall that $W_t=\rho \cdot  B_t+\bar\rho\cdot  \tilde B_t$. If   there is a constant $ J$ and a process $\theta_t=\theta(t,Y_t)$  such that  $\theta^2_t\cdot \Lambda_T<\infty$ and 
\beq
\label{Hobson1}
J+\theta\cdot W_T+\Big(\rho \theta\lambda+\frac 12\bar \rho^2\theta^2-\frac 12 \lambda^2\Big)\cdot \Lambda_T=0,  
\eeq
then $ J$, $\xi_t:= \bar \rho_t\theta_t$, and $H_t:= \rho_t \theta_t$ form a solution of the equation (\ref{Hequation}). 
  
Applying the It\^o formula with  $f\in C^{1,2}$ solving ({\ref{PDE}}) and $Y$ given by (\ref{2+}) we get that 
\bean
0&=&f(0,Y_0)+\beta f' \cdot W_T+\Big (\dot f + \alpha f'+\frac 12 \beta^2f''\Big)\cdot \Lambda_T\\
&=&
f(0,Y_0)+\beta f'  \cdot W_T+\Big (\rho\lambda\beta f' +\frac 12 \bar \rho^2\beta^2f'^2 -\frac 12 \lambda^2\Big)\cdot \Lambda_T.
\eean
Thus, (\ref{Hobson1}) holds with $ J=f(0,Y_0)$ and $\theta_t=\beta (t,Y_t) f'(t,Y_t)$. \fdem

\section{Minimal entropy EMM for the model (\ref{1}) -- (\ref{3})}

\begin{theorem} The density of minimal entropy EMM in the model (\ref{1}) -- (\ref{3}) has the form 
\beq
\label{Zo}
Z^o_T=\frac{dQ^o}{dP} = \exp\left\{-\int_0^T \frac{\mu + \kappa t}{\sigma} dB_t - \frac 12\int_0^T\frac{(\mu + \kappa t)^2}{\sigma^2} dt\right\}.
\eeq
Moreover, 
with respect to the measure $Q^o=Z^o_T P$ this  model admits the representation 
\bea
\label{11}
dS_t&=&S_tY_t dB^o_t,\\
\label{21}
dY_t&=&Y_y\Big (\theta  + \frac{1}{2}\beta^2+\alpha \ln \sigma  - \beta\rho\frac{\mu + \kappa t}{\sigma} - \alpha \ln Y_t\Big)dt+ \beta Y_t dW^o_t, 
\eea 
where $\E^o[ B^o_tW^o_t]=\rho t$.
\end{theorem}
{\sl Proof.} Let us denote by $Y_t:=\sigma e^{V_t}$. By the It\^ o formula 
$$
dY_t = Y_t\Big(\theta + \frac{1}{2}\beta^2 +\alpha \ln \sigma  - \alpha\ln Y_t\Big)dt + \beta Y_tdW_t.
$$
In the notations of the  previous section we have:  
$$
\lambda (t) = \frac{\mu + \kappa t}{\sigma},  \quad  \alpha(y) = y\left(\theta + \frac{1}{2}\beta^2 + \alpha \ln \sigma\right) - \alpha y\ln y, \quad   \beta (y)=\beta y. 
$$ 
It is easily seen that in our specific case  that the  solution of  the terminal problem (\ref{PDE})  does not depend on $y$. That is, it satisfies the terminal value problem for the ordinary differential equation 
\beq
\label{PDE1}
\dot f+\frac 12\lambda^2=0, \quad f(T)=0, 
\eeq
and admits the explicit expression 
$$
f(t)=\frac{(\mu + \kappa (T-t))^3-\mu^3}{6\sigma^2\kappa}. 
$$
Also, $H=0$ and $\xi=0$. The density of the minimal entropy EMM is given by (\ref{Zo}) and $J=f(0)$.

 \section{Simulations}
We tested the studied model on real market data from the site https: //finance.yahoo.com/  using the daily closing prices of shares of the Alphabet  (ticker: GOOG) from Dec. 15, 2022 to Apr. 21, 2023 and evaluating the  American call option GOOG230421C00090000, the market prices of which  were taken from 
 https://www.barchart.com/options/. The expiration date $T$   was Apr. 21, 2023 and the   strike $K=90$.
 \begin{figure}[htb]
    \centering
    \includegraphics[width=0.9\textwidth]{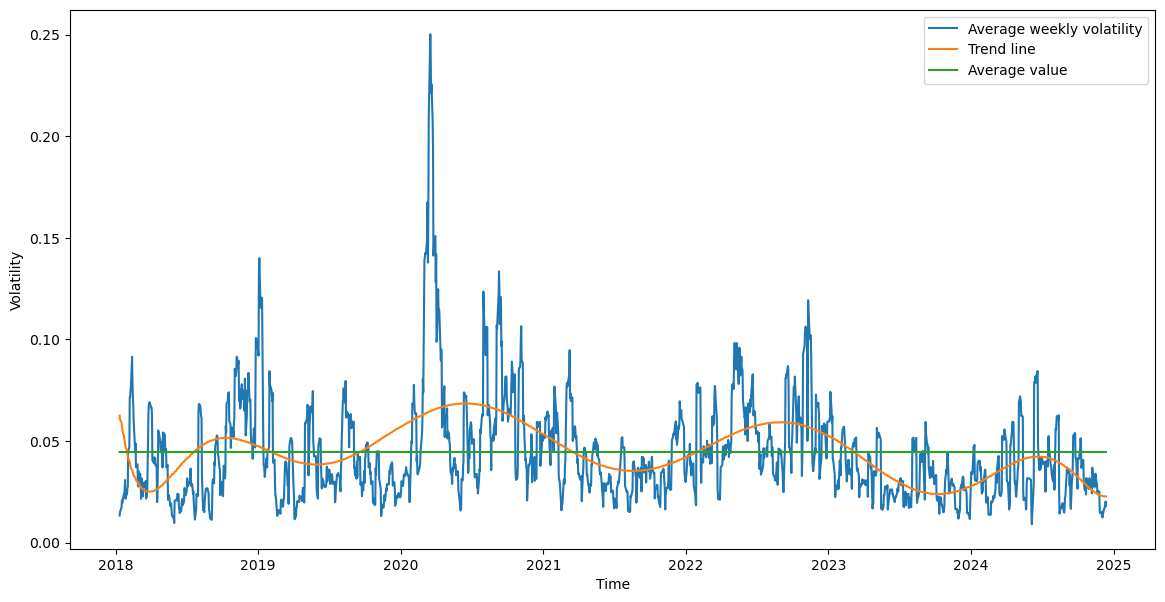}
    \caption{\small  The average weekly volatility  for AAPL stock. The share price data are taken from the website {https://finance.yahoo.com/} in the period from Jan. 1, 2018 to Dec. 1, 2024.}
    \label{fig:my_label}
\end{figure}

We use the Monte Carlo method with the control variate as described in \cite{Glasserman}, Ch. 4. We applied the Euler scheme to simulate the trajectories of (\ref{11}) and (\ref{21}) and used the following approximations:
$$
\begin{cases}
\hat{Y}_{t_{i+1}} = \hat{Y}_{t_i} + \hat{Y}_{t_i}(\theta + \frac{1}{2}\beta^2 + \alpha \ln \sigma - \beta\rho\frac{\mu + \kappa t_{i}}{\sigma} - \alpha\ln \hat{Y}_{t_i}) \Delta t_{i+1} + \beta\hat{Y}_{t_i} \Delta W^o_{t_{i+1}}, \\[10pt]
\hat{X}_{t_{i+1}} = \hat{X}_{t_i} - \frac{1}{2}\hat{Y_{t_{i}}}^2\Delta t_{i+1} + \hat{Y}_{t_i}\Delta B^o_{t_{i+1}},  \\[10pt]
\hat{S}_{t_{i+1}} = \exp\{\hat{X}_{t_{i+1}}\},
\end{cases}
$$
where $0=t_0 < t_1 < ... < t_m$ is the time grid, $\Delta t_{i+1} := t_{i+1} - t_{i}$, the increment $\Delta W^o_{t_{i+1}} := W^o_{t_{i+1}} - W^o_{t_i}$, $W^o_t: = \rho B^o_t  + \bar \rho \tilde B^o_t$, $B^o_t$ and $\tilde B^o_t$ are independent Wiener processes. 

In Fig. 8.2  we represent  the result of our model after selecting  parameters. The relative error is 0.117, the coefficient of determination $R^2$ = 0.89 indicating a relatively good quality of the model.

\begin{figure}
    \centering
    \includegraphics[width=0.95\textwidth]{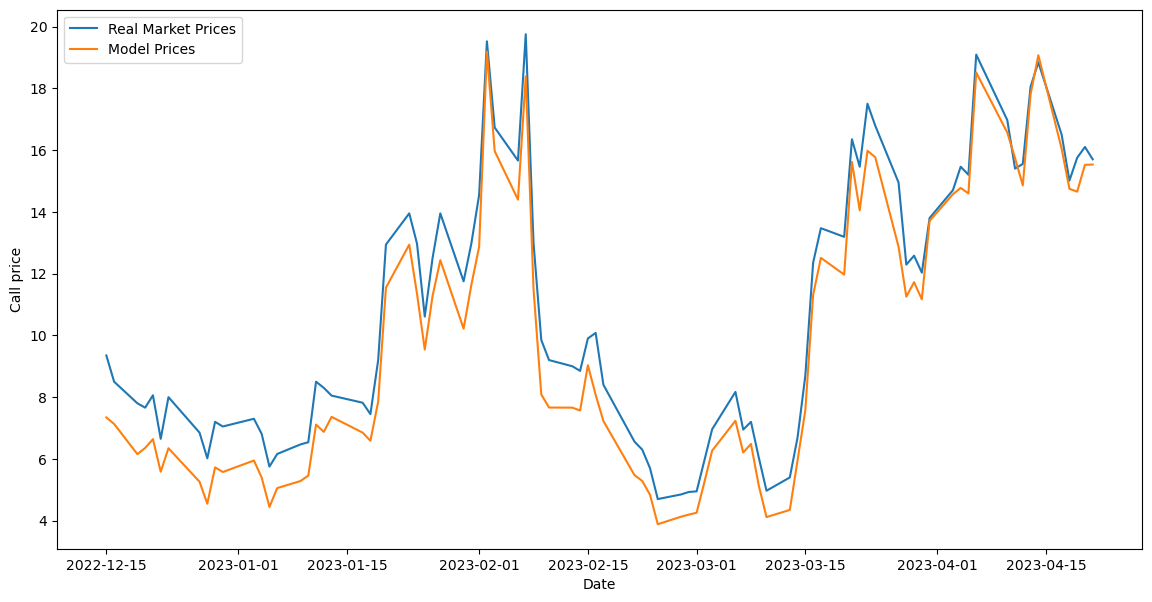}
    \caption{\small Plot of American call option prices with strike $K = 90$ and expiry date $T = 21.04.2023$ calculated by the model and the realised market prices of this instrument. The model parameters are:  $Y_0 = 0.3, \mu = 0.25, \kappa = 0, \sigma = 0.5, \alpha = 0.75, \theta = 0.1, \beta = 0.2, \rho = -0.5$.}
    \label{fig:my_label}
\end{figure}

\section{Conclusion}
In this note we present a short theoretical prerequisites for a search of the minimal entropy equivalent martingale using the Hobson construction. As an application, we 
 obtained an explicit expression for such a measure  for  a specific stochastic volatility model referred in the literature as the exponential Ornstein--Uhlenbeck model in the form suitable for pricing contingent claims.

\section*{Competing interests}
The authors declare no competing interests.

%

\newpage

\end{document}